\documentclass{amsart}
\usepackage{graphicx}
\usepackage{amscd}
\usepackage{amsmath}
\usepackage{amsxtra}
\usepackage{amsfonts}
\usepackage{amssymb}

\newtheorem{theorem}{Theorem}[section]

\newtheorem{lemma}[theorem]{Lemma}
\newtheorem{proposition}[theorem]{Proposition}
\newtheorem{conjecture}[theorem]{Conjecture}
\theoremstyle{definition}
\newtheorem{definition}[theorem]{Definition}
\newtheorem{remark}[theorem]{Remark}
\newtheorem*{problem}{Problem}
\newtheorem{example}[theorem]{Example}
\theoremstyle{remark}

\renewcommand{\theclaim}{\textup{\theclaim}}

\newtheorem*{acknowledgements}{Acknowledgements}

\numberwithin{equation}{section}

\def\openone

{\mathchoice

{\hbox{\upshape \small1\kern-3.3pt\normalsize1}}

{\hbox{\upshape \small1\kern-3.3pt\normalsize1}}

{\hbox{\upshape \tiny1\kern-2.3pt\SMALL1}}

{\hbox{\upshape \Tiny1\kern-2pt\tiny1}}}

\makeatletter

\newbox\ipbox

\newcommand{\ip}[2]{\left\langle #1\,|\,#2\right\rangle}
\newcommand{\vect}[2]{\left[\begin{array}{cc}
#1\\ #2 \end{array}\right]}
\newcommand{\diracb}[1]{\left\langle #1\mathrel{\mathchoice

{\setbox\ipbox=\hbox{$\displaystyle \left\langle\mathstrut
#1\right.$}

\vrule height\ht\ipbox width0.25pt depth\dp\ipbox}

{\setbox\ipbox=\hbox{$\textstyle \left\langle\mathstrut
#1\right.$}

\vrule height\ht\ipbox width0.25pt depth\dp\ipbox}

{\setbox\ipbox=\hbox{$\scriptstyle \left\langle\mathstrut
#1\right.$}

\vrule height\ht\ipbox width0.25pt depth\dp\ipbox}

{\setbox\ipbox=\hbox{$\scriptscriptstyle \left\langle\mathstrut
#1\right.$}

\vrule height\ht\ipbox width0.25pt depth\dp\ipbox}

}\right. }

\newcommand{\dirack}[1]{\left. \mathrel{\mathchoice

{\setbox\ipbox=\hbox{$\displaystyle \left.\mathstrut
#1\right\rangle$}

\vrule height\ht\ipbox width0.25pt depth\dp\ipbox}

{\setbox\ipbox=\hbox{$\textstyle \left.\mathstrut
#1\right\rangle$}

\vrule height\ht\ipbox width0.25pt depth\dp\ipbox}

{\setbox\ipbox=\hbox{$\scriptstyle \left.\mathstrut
#1\right\rangle$}

\vrule height\ht\ipbox width0.25pt depth\dp\ipbox}

{\setbox\ipbox=\hbox{$\scriptscriptstyle \left.\mathstrut
#1\right\rangle$}

\vrule height\ht\ipbox width0.25pt depth\dp\ipbox}

} #1\right\rangle}

\newcommand{\cj}[1]{\overline{#1}}

\newcommand{\bz}{\mathbb{Z}}
\newcommand{\br}{\mathbb{R}}
\newcommand{\bc}{\mathbb{C}}

\newcommand{\bn}{\mathbb{N}}

\usepackage{graphicx}

\def\blfootnote{\xdef\@thefnmark{}\@footnotetext}

\newcommand{\wdots}{\ldots}
\newcommand{\proj}{\operatorname{proj}}

\begin{document}
\title[Fourier frequencies in affine IFS]{Fourier frequencies in affine iterated function systems}
\author{Dorin Ervin Dutkay}
\blfootnote{Research supported in part by the National Science
Foundation DMS 0457491}
\address[Dorin Ervin Dutkay]{\protect\raggedright
Department of Mathematics,
University of Central Florida,
4000 Central Florida Blvd., P.O. Box 161364,
	Orlando, FL 32816-1364}
\email{ddutkay@mail.ucf.edu}
\author{Palle E.T. Jorgensen}
\address[Palle E.T. Jorgensen]{Department of Mathematics, 
The University of Iowa, 14 MacLean Hall, Iowa City, IA 52242-1419, USA}
\email{jorgen@math.uiowa.edu}
\begin{abstract}
 We examine two questions regarding Fourier frequencies for a class of iterated function systems (IFS). These are iteration limits arising from a fixed finite families of affine and contractive mappings in $\br^d$, and the ``IFS'' refers to such a finite system of transformations, or functions. The iteration limits are pairs $(X, \mu)$ where $X$ is a compact subset of $\br^d$, (the support of $\mu$) and the measure $\mu$ is a probability measure determined uniquely by the initial IFS mappings, and a certain strong invariance axiom. The two questions we study are: (1) existence of an orthogonal Fourier basis in the Hilbert space $L^2(X,\mu)$; and (2) explicit constructions of Fourier bases from the given data defining the IFS.
\end{abstract}

\maketitle \tableofcontents

\section{Introduction}\label{introduction}
\par
Motivated in part by questions from wavelet theory, there has been a set of recent advances in a class of spectral problems from iterated function systems (IFS) of affine type. The geometric side of an IFS is a pair $(X, \mu)$ where $X$ is a compact subset of $\br^d$, (the support of $\mu$) and the measure $\mu$ is a probability measure determined uniquely by the initial IFS mappings, and a certain strong invariance property. In this paper, we examine two questions regarding Fourier frequencies for these iterated function systems (IFS): (1) When do we have existence of an orthogonal Fourier basis in the Hilbert space $L^2(X, \mu)$; and, when we do, (2) explicitly, what are the Fourier frequencies of these orthonormal bases in terms of the data that defines the iterated function system? Our main result, Theorem \ref{thmain},	 shows that existence in (1) follows from geometric assumptions that are easy to check, and it is a significant improvement on earlier results in the literature.  Our approach uses a new idea from dynamics, and it allows us to also answer (2).

\par
By a Fourier basis in $L^2(X,\mu)$ we mean a subset $\Lambda$ of $\br^d$ such that the functions $\{e_\lambda\, |\, \lambda \in \Lambda\}$ form an orthogonal basis in $L^2(\mu)$. Here $e_\lambda (x) := \exp( 2 \pi i \lambda \cdot x)$. The functions $e_\lambda$ are restricted from $\br^d$ to $X$. (The factor $2 \pi$ in the exponent is introduced for normalization purposes only.)
\par
So far Fourier bases have been used only in the familiar and classical context of compact abelian groups; see, e.g., \cite{Kat04}. There, as is well known, applications abound, and hence it is natural to attempt to extend the fundamental duality principle of Fourier bases to a wider category of sets $X$ which are not groups and which in fact carry much less structure. Here we focus on a particular such class of subsets $X$ in $\br^d$ which are IFS attractors. Our present paper focuses on the theoretical aspects which we feel are of independent interest, but we also allude to applications. 
\par
Since $X$ and its boundary are typically fractals in the sense of \cite{Man04}, their geometry and structure do not lend themselves in an obvious way to Fourier analysis. (Recall \cite{Man04} that some fractals model chaos.) To begin with, the same set $X$ may arise in more than one way as a limiting object. It will be known typically from some constructive algorithm. While each finite algorithmic step can readily be pictured, not so for the iteration limit!  And from the outset it may not even be clear whether or not a particular $X$ is the attractor of an iterated function system (IFS); see, e.g., \cite{LaFr03, Fal03, Jor06,Bea65, BCMG04}. Moreover, far from all fractals fall in the affine IFS class. But even the affine class of IFSs has a rich structure which is not yet especially well understood.
\par
The presence of an IFS structure for some particular set $X$ at least implies a preferred self-similarity; i.e., smaller parts of $X$ are similar to its larger scaled parts, and this similarity will be defined by the maps from the IFS in question. When $X$ is the attractor of a given contractive IFS $(\tau_i)$, then by \cite{Hut81}, there is a canonical positive and strongly invariant measure $\mu$ which supports $X$. But even in this case, a further difficulty arises, addressed in Section \ref{sect4} below. 

\par
As illustrated with examples in Section \ref{sect5} below, the geometric patterns for a particular $X$ might not at all be immediately transparent. For a given $X$, the problem is to detect significant patterns such as self-similarity, or other ``hidden structures'' (see, e.g., \cite{CuSm02, Sma05}); and Fourier frequences, if they can be found, serve this purpose. In addition, if $X$ does admit a Fourier basis, this allows us to study its geometry and its symmetries from the associated spectral data. In that case, standard techniques from Fourier series help us to detect ``hidden'' structures and patterns in $X$. 
\par
However, we caution the reader that recent work of Strichartz \cite{Str05} shows that a number of ``standard'' results from classical Fourier series take a different form in the fractal case.
\par
In the next section we give definitions and recall the basics from the theory of iteration limits; i.e., metric limits which arise from a fixed finite family of affine and contractive mappings in $\br^d$, and the ``IFS'' refers to such a finite system of transformations.
\par
There are a number of earlier papers \cite{JoPe98,DuJo05,LaWa02,Str00,LaWa06} which describe various classes of affine IFSs $(X, \mu)$ for which an orthogonal Fourier basis exists in $L^2(X, \mu)$. It is also known \cite{JoPe98} that if the affine IFS $(X, \mu)$ is the usual middle-third Cantor set, then no such Fourier basis exists; in fact, in that case there can be no more than two orthogonal Fourier exponentials $e_\lambda$ in $L^2(X, \mu)$. Nonetheless, the present known conditions which imply the existence of an orthogonal Fourier basis have come in two classes, an algebraic one (Definition \ref{defhada} below) and an analytic assumption. Our main result, Theorem \ref{thmain}, shows that the analytic condition can be significantly improved. We also conjecture that the algebraic condition is sufficient (see Conjecture \ref{conj}).

\section{Definitions and preliminaries}\label{defandpre}

    The definitions below serve to make precise key notions which we need
to prove the main result (Theorem \ref{thmain}). In fact they are needed in relating
the intrinsic geometric features of a given affine IFS $(X, \mu)$ to the
spectral data for the corresponding Hilbert space $L^2(X, \mu)$.  Our paper
focuses on a class of affine IFSs which satisfies a certain symmetry
condition (Definition \ref{defhada}). This condition involves a pair of IFSs in
duality, and a certain complex Hadamard matrix. While these duality systems
do form a restricted class, their study is motivated naturally by our
recursive approach to building up a Fourier duality. Moreover, our recursive
approach further suggests a certain random-walk model which is built
directly on the initial IFS. We then introduce a crucial notion of invariant
sets for this random walk (Definition \ref{Def2.9}). The corresponding transition
probabilities of the random walk are defined in terms of the Hadamard matrix
in Definition \ref{defhada}, and it lets us introduce a discrete harmonic analysis, a
Perron--Frobenius operator and associated harmonic functions (Definition
\ref{Def2.6}).  The interplay between these functions and the invariant sets is made
precise in Propositions \ref{propfinite} and \ref{propsupp}, and Theorem \ref{thcora}.

\begin{definition}
\label{Def2.1}
A probability measure $\mu$ on $\br^d$ is called a {\it spectral measure} if there exists a subset $\Lambda$ of $\br^d$ such that the family of exponential functions $\{e^{2\pi i\lambda\cdot x}\,|\, \lambda\in\Lambda\}$ is an orthonormal basis for $L^2(\mu)$. In this case, the set $\Lambda$ is called a {\it spectrum} of the measure $\mu$.
\end{definition}

       It was noted recently in \cite{LaWa06} that the axiom which defines
spectral measures $\mu$ implies a number of structural properties for $\mu$,
as well as for the corresponding spectrum $\Lambda = \Lambda(\mu)$: e.g.,
properties regarding discreteness and asymptotic densities for $\mu$, and
intrinsic algebraic relations on the configuration of vectors in $\Lambda$.

       Our present paper deals with the subclass of spectral measures that
can arise from affine IFSs.

\begin{definition}
Let $Y$ be a complete metric space. Following \cite{Hut81} we say that a finite family $(\tau_i)_{i=1,N}$ of contractive mappings in $Y$ is an iterated function system (IFS). Introducing the Hausdorff metric on the set of compact subsets $K$ of $Y$, we get a second complete metric space, and we note that the induced mapping
$$K\mapsto\bigcup_{i=1}^N\tau_i(K),$$ 
is contractive. By Banach's theorem, this mapping has a unique fixed point, which we denote $X$; and we call $X$ the {\it attractor} for the IFS. It is immediate by restriction that the individual mappings $\tau_i$ induce endomorphisms in $X$, and we shall denote these restricted mappings also by $\tau_i$.
\par
For IFSs where the mappings $\tau_i$ are affine as in (\ref{eqtaub}) below, we talk of {\it affine} IFSs. In this case, the ambient space is $\br^d$. 
\end{definition}

Let $R$ be a $d\times d$ expansive integer matrix, i.e., all entries are integers and all eigenvalues have absolute value strictly bigger than one. For a point $b\in\bz^d$ we define the function 
\begin{equation}\label{eqtaub}
\tau_b(x):=R^{-1}(x+b)\quad(x\in\br^d).
\end{equation}
\par
For a finite subset $B\subset\bz^d$ we will consider the iterated function system $(\tau_b)_{b\in B}$. We denote by $N$ the cardinality of $ B$. We will assume also that $0\in B$.
\par
The fact that the matrix $R$ is expansive implies that there exists a norm on $\br^d$ for which the maps $\tau_b$ are contractions.
\par
There exist then a unique compact set $X_B$, called the {\it attractor} of the IFS, with the property that 
$$X_B=\bigcup_{b\in B}\tau_b(X_B).$$
Moreover, we have the following representation of the attractor: 
$$X_B=\left\{\sum_{k=1}^\infty R^{-k}b_k\biggm| b_k\in B\mbox{ for all }k\geq1\right\}.$$
\par
There exists a unique invariant probability measure $\mu_B$ for this IFS, i.e., for all bounded continuous functions on $\br^d$,
\begin{equation}
\int f\,d\mu_B=\frac{1}{N}\sum_{b\in B}\int f\circ\tau_b\,d\mu_B.
\label{eqmuB}
\end{equation}
Moreover, the measure $\mu_B$ is supported on the attractor $X_B$. We refer to \cite{Hut81} for details.
\par
Following earlier results from \cite{JoPe98,Str00,LaWa02,DuJo05,LaWa06}, in order to obtain Fourier bases for the measure $\mu_B$, we will impose the following algebraic condition on the  pair $(R,B)$:
\begin{definition}\label{defhada}
Let $R$ be a $d\times d$ integer matrix, $B\subset\bz^d$ and $L\subset\bz^d$ having the same cardinality as $B$, $\#B=\#L=:N$. We call $(R,B,L)$ a {\it Hadamard triple} if 
the matrix
$$\frac{1}{\sqrt{N}}(e^{2\pi iR^{-1}b\cdot l})_{b\in B,l\in L}$$
is unitary.
\end{definition}

\par
We will assume throughout the paper that $(R,B,L)$ is a Hadamard triple. 

\begin{remark}\label{reminco}
Note that if $(R,B,L)$ is a Hadamard triple, then no two elements in $B$ are congruent modulo $R\bz^d$, and no two elements in $L$ are congruent modulo $R^T\bz^d$. 
\par
Indeed, if $b,b'\in B$ satisfy $b-b'=Rk$ for some $k\in\bz^d$ then, since $L\subset\bz^d$,
$$e^{2\pi i R^{-1}b\cdot l}=e^{2\pi iR^{-1}b'\cdot l}\quad(l\in L),$$
so the rows $b$ and $b'$ of the matrix in Definition \ref{defhada} cannot be orthogonal.
\end{remark}
\par
We conjecture that the existence of a set $L$ such that $(R,B,L)$ is a Hadamard triple is sufficient to obtain orthonormal bases of exponentials in $L^2(\mu_B)$. 
\begin{conjecture}\label{conj}
Let $R$ be a $d\times d$ expansive integer matrix, $B$ a subset of $\bz^d$ with $0\in B$. Let $\mu_B$ be the invariant measure of the associated IFS $(\tau_b)_{b\in B}$. If there exists a subset $L$ of $\bz^d$ such that $(R,B,L)$ is a Hadamard triple and $0\in L$ then $\mu_B$ is a spectral measure. 
\end{conjecture}
\par
We will prove in Theorem \ref{thmain} that the conjecture is true under some extra analytical assumptions, thus extending the known results from \cite{JoPe98,Str00,LaWa02,DuJo05,LaWa06}.

\subsection{Path measures}\label{pathmeasures}
To analyze the measure $\mu_B$ we will use certain random-walk (or ``path'') measures $P_x$ which are directly related to the Fourier transform $\hat\mu_B$ of the invariant measure. Most of the results in Sections \ref{pathmeasures} and \ref{invariantsets} are essentially contained in \cite{CoRa90,CCR96,DuJo05}. We include them here for the convenience of the reader.
\par
Define the function 
$$W_B(x)=\left|\frac{1}{N}\sum_{b\in B}e^{2\pi ib\cdot x}\right|^2\quad(x\in\br^d).$$
This function appears if one considers the Fourier transform of equation (\ref{eqmuB}):
$$\left|\hat\mu_B(x)\right|^2=W_B\left((R^T)^{-1}x\right)\left|\hat\mu_B\left((R^T)^{-1}x\right)\right|^2,\quad(x\in\br^d).$$
\par
The elements of $L$ and the transpose $S:=R^T$ will define another iterated function system 
$$\tau_l(x)=S^{-1}(x+l)\quad(x\in\br^d,l\in L).$$
\par
We underline here that we are interested in the measure $\mu_B$ associated to the iterated function system $(\tau_b)_{b\in B}$, and the main question is whether this is a spectral measure. The iterated function system $(\tau_l)_{l\in L}$ will only help us in constructing the basis of exponentials. 
\par
The unitarity of the matrix in Definition \ref{defhada} implies (see \cite{LaWa02}, \cite{DuJo05}) that
\begin{equation}\label{eqqmf}
\sum_{l\in L}W_B(\tau_lx)=1\quad(x\in\br^d).
\end{equation}

\begin{remark}
\label{RemarkApr25.pound}
     The reader will notice that in our analysis of the iteration steps, our
measure $\mu_B$ in \textup{(\ref{eqmuB})} is chosen in such a way that each of the branches
in the iterations is given equal weight $1/N$. There are a number of reasons
for this.

But first recall the following known theorem from \cite{Hut81} to the effect that
for every IFS $(\tau_b)_{b\in B}$, $b \in B$, $N = \# B$, and for every $N$-configuration of
numerical weights $(p_b)_{b\in B}$,  $p_b > 0$, with $\sum_{b\in B} p_b = 1$, there is a unique
$(p_b)$-distributed probability measure $\mu_{p,B}$ with support $X_B$. This measure
$\mu_{p,B}$ is determined uniquely by the equation
\[
  \mu_{p,B} = \sum_{b\in b}  p_b \mu_{p,B}\circ \tau_b^{-1} .
\]
Since our focus is on spectral measures \textup{(}Definition \textup{\ref{Def2.1})}, it is natural to
restrict attention to the case of equal weights, i.e., to $p_b = 1/N$.
\end{remark}

Another reason for this choice is a conjecture by \L aba and Wang \cite{LaWa02}, as
well as the following lemma.

\begin{lemma}
\label{LemmaApr25.pound}
Set  $$W_{p,B}(x):=\left|\sum_{b\in B} p_b e^{2\pi i b\cdot x}\right|^2,$$ and assume that
$$\sum_{l\in L} W_{p,B}(\tau_l(x))=1$$ for some dual IFS 
$$\tau_l(x)=(R^T)^{-1}(x+l),\quad (x\in\br^d,l \in L),$$ with $\#L=N$. Then $p_b=1/N$ for all $b\in B$. 
\end{lemma}
 
\begin{proof}
Expanding the modulus square and changing the order of sumation, we get that for all $x\in\br^d$,
$$\sum_{b,b'\in B}p_bp_{b'}e^{2\pi i R^{-1}(b-b')\cdot x}\sum_{l\in L}e^{2\pi iR^{-1}(b-b')\cdot l}=1$$
The constant term on the left must be equal to $1$, so
$$\sum_{b\in B}Np_b^2=1.$$ 
Since $\sum_{b\in B}p_b=1$, this will imply that we have equality in a Schwarz inequality, so $p_b=1/N$ for all $b\in B$.  
\end{proof}

The relation (\ref{eqqmf}) can be interpreted in probabilistic terms: $W_B(\tau_lx)$ is the probability of transition from $x$ to $\tau_lx$. This interpretation will help us define the path measures $P_x$ in what follows.
\par

Let $\Omega=\{(l_1l_2\dots)\,|\,l_n\in L\mbox{ for all }n\in\bn\}=L^{\bn}$. Let $\mathcal{F}_n$ be the sigma-algebra generated by the cylinders depending only on the first $n$ coordinates.
\par
There is a standard way due to Kolmogorov of using the system $(\br^d, (\tau_b)_{b\in B})$ to generate a path space $\Omega$, and an associated family of path-space measures $P_x$, indexed by $x\in\br^d$. Specifically, using the weight function $W_B$ in assigning conditional probabilities to random-walk paths, we get for each $x\in\br^d$ a Borel measure $P_x$ on the space of paths originating in $x$. For each $x$, we consider paths originating at $x$, and governed by the given IFS. The transition probabilities are prescribed by $W_B$; and passing to infinite paths, we get the measure $P_x$. We shall refer to this $(P_x)_{x\in\br^d}$ simply as the path-space measure, or the {\it path measure} for short.
\par
For each $x\in\br^d$ we can define the measures
$P_x$ on $\Omega$ as follows. For a function $f$ on $\Omega$ which depends only on the first $n$ coordinates
$$\int_{\Omega}f\,dP_x=\sum_{\omega_1,\dots ,\omega_n\in L}W_B(\tau_{\omega_1}x)W_B(\tau_{\omega_2}\tau_{\omega_1}x)\cdots W_B(\tau_{\omega_n}\cdots \tau_{\omega_1}x)f(\omega_1,\dots ,\omega_n).$$
In particular, when the first $n$ components are fixed $l_1,\dots ,l_n\in L$,
\begin{equation}\label{eqcyl}
P_x(\{(\omega_1\omega_2\wdots )\in\Omega\,|\,\omega_1=l_1,\dots ,\omega_n=l_n\})=\prod_{k=1}^nW_B(\tau_{l_k}\cdots \tau_{l_1}x).
\end{equation}
\par
Define the transfer operator 
$$R_Wf(x)=\sum_{l\in L} W_B(\tau_lx)f(\tau_lx)\quad(x\in\br^d).$$

\begin{definition}
\label{Def2.6}
A measurable function $h$ on $\br^d$ is said to be $R_W${\it -harmonic} if  $R_W h = h$. A measurable function $V$ on $\br^d\times\Omega$ is said to be a cocycle if it satisfies the following covariance property:
\begin{equation}\label{eqinv}
V(x,\omega_1\omega_2\wdots )=V(\tau_{\omega_1}x,\omega_2\omega_3\wdots )\quad(\omega_1\omega_2\wdots \in\Omega).
\end{equation}
\end{definition}
\par
In the following we give a formula for all the bounded $R_W$-harmonic functions. The result expresses the bounded $R_W$-harmonic functions in terms of a certain boundary integrals of cocycles, and it may be viewed as a version of the Fatou--Markoff--Primalov theorem.
\par
If $h$ is a bounded measurable $R_W$-harmonic function on $\br^d$, then, for all $x\in\br^d$, the functions 
$$(\omega_1,\dots ,\omega_n)\mapsto h(\tau_{\omega_n}\cdots \tau_{\omega_1}x)$$
define a bounded martingale. By Doob's martingale theorem, one obtains that the following limit exists $P_x$-a.e.:
\begin{equation}\label{eqconv}
\lim_{n\rightarrow\infty}h(\tau_{\omega_n}\cdots \tau_{\omega_1}x)=:V(x,\omega),\mbox{ for }P_x\mbox{-a.e.\ }\omega\in\Omega,
\end{equation} where $V(x,\cdot):\omega\rightarrow\bc$ is some bounded function on $\Omega$. Moreover, $V$ is a cocycle.
\par
We formalize this conclusion in a lemma.
\begin{lemma}\label{lemhaco}
If $h$ is a bounded $R_W$-harmonic function, then the associated function $V$ from \textup{(\ref{eqconv})} is well defined, it is bounded and measurable; and it is a cocycle. Conversely, if $V:\br^d\times\Omega\rightarrow\bc$ is a bounded measurable function satisfying \textup{(\ref{eqinv})}, then the function 
\begin{equation}\label{eqh_V}
h_V(x):=P_x(V(x,\cdot))\quad(x\in\br^d),
\end{equation}
defines a bounded function on $\br^d$ such that $R_Wh_V=h_V$, and such that relation \textup{(\ref{eqconv})} is satisfied with $h=h_V$. 
\end{lemma}
\par
Next we show that the family of measures $x\mapsto P_x$ is weakly continuous. More precisely, we have the following result.

\begin{proposition}\label{propcontpx}\cite[Proposition 5.2]{CoRa90}
Let $U$ be a bounded measurable function on $\Omega$. Then there exists a constant $0\leq D<\infty$ such that
$$|P_x(U)-P_y(U)|\leq D|x-y|\|U\|_\infty\quad(x,y\in\br^d).$$
\end{proposition}
While the main ideas are contained in \cite{CoRa90}, we include the proof for the benefit of the reader;
our version covers affine matrix operations for contraction, extending
the one-dimensional dyadic case in \cite{CoRa90}.
\begin{proof} 
Let $x,y\in\br^d$. For $\omega_1\wdots \omega_n\in L^n$ and $1\leq p\leq n$, define $W_{\omega,p}(x):=W_B(\tau_{\omega_p}\cdots \tau_{\omega_1}x)$, and
$$\delta_n(x,y):=\sum_{\omega_1\wdots \omega_n\in L^n}|W_{\omega,n}(x)\cdots W_{\omega,1}(x)-W_{\omega,n}(y)\cdots W_{\omega,1}(y)|.$$
We have, using equation (\ref{eqqmf}),
$$\delta_n(x,y)\leq\sum_{\omega_1\wdots \omega_n\in L^n}|W_{\omega,n}(x)-W_{\omega,n}(y)|W_{\omega,n-1}(x)\cdots W_{\omega,1}(x)+\delta_{n-1}(x,y)$$
$$\leq Mc^n|x-y|+\delta_{n-1}(x,y),$$
where $c$ is the contraction constant for the maps $\tau_{l}$, $l\in L$, and $M$ is a Lipschitz constant for $W_B$.
\par
{}From this we obtain
$$\delta_n(x,y)\leq M|x-y|\sum_{k\geq 1}c^k.$$
This proves the result in the case when $U$ depends only on a finite number of coordinates.
\par
In the general case, define $Q:=\frac{1}{2}(P_x+P_y)$, and let $U_n$ be the conditional expectation $\mathbb{E}_Q[U|\mathcal{F}_n]$. The functions $U_n$, $n\geq1$, are bounded by $\|U\|_\infty$ and the sequence converges 
$Q$-a.e., and so $P_x$ and $P_y$-a.e., to $U$. It follows from the previous estimate that
$$|P_x(U_n)-P_y(U_n)|\leq\|U\|_\infty\delta_n(x,y)\leq D|x-y|\|U\|_\infty.$$
The result is obtained by applying Lebesgue's dominated convergence theorem.
\end{proof}

\subsection{Invariant sets}\label{invariantsets}
In the following, we will work with the affine system $(\tau_l)_{l\in L}$, and with the weight function $W_B$. Given this pair, we introduce a notion of invariant sets as introduced in \cite{CoRa90,CCR96,CHR97}. We emphasize that ``invariance'' depends crucially on the chosen pair. The reason for the name ``invariance'' is that the given affine system and the function $W_B$ naturally induce an associated  random walk on points in $\br^d$ as described before. 
\par
Let $x$ and $y$ be points in $\br_d$ and suppose $y =\tau_l (x)$ for some $l\in L$. We then say that $W_B(y)$ represents the probability of a transition from $x$ to $y$. Continuing this with paths of points, we then arrive at a random-walk model, and associated trajectories, or paths. An orbit of a point $x$  consists of the closure of the union of those trajectories beginning at $x$ that have positive transition probability between successive points. A closed set $F$ will be said to be invariant if it contains all its orbits starting in $F$. Note in particular that every (closed) orbit is an invariant set.
\par
We now spell out these intuitive notions in precise definitions.
\begin{definition}
\label{Def2.9}
For $x\in\br^d$, we call a {\it trajectory} of $x$ a set of points $$\{\tau_{\omega_n}\cdots \tau_{\omega_1}x\,|\,n\geq1\}$$ where 
$\{\omega_n\}_n$ is a sequence of elements in $L$ such that $W_B(\tau_{\omega_n}\cdots \tau_{\omega_1}x)\neq0$ for all $n\geq1$. We denote by $\mathcal{O}(x)$ the union of all trajectories of $x$ and the closure $\cj{\mathcal{O}(x)}$ is called the {\it orbit} of $x$. If $W_B(\tau_lx)\neq0$ for some $l\in L$ we say that the {\it transition} from $x$ to $\tau_lx$ is possible.
\par
A closed subset $F\subset\br^d$ is called {\it invariant} if it contains the orbit of all of its points. An invariant subset is called {\it minimal} if it does not contain any proper invariant subsets.
\end{definition}
A closed subset $F$ is invariant if, for all $x\in F$ and $l\in L$ such that $W_B(\tau_lx)\neq 0$, it follows that $\tau_lx\in F$.
\par
Since the orbit of any point is an invariant set,
a closed subset $F$ is minimal if and only if $F=\cj{\mathcal{O}(x)}$ for all $x\in F$. By Zorn's lemma, every invariant subset contains a minimal subset.
\begin{proposition}\label{propmini}
If $F_1$ is a closed invariant subset and $F_2$ is a compact minimal invariant subset of $\br^d$ then either $F_1\cap F_2=\emptyset$ or $F_2\subset F_1$.
\end{proposition}

\begin{proof}
Indeed, if $x\in F_1\cap F_2$ then $F_2=\cj{\mathcal{O}(x)}\subset F_1$.
\end{proof}
\par
\begin{proposition}\label{propnf}
Let $F$ be a compact invariant subset. Define 
$$N(F):=\{\omega\in\Omega\,|\, \lim_{n\rightarrow\infty}d(\tau_{\omega_n}\cdots \tau_{\omega_1}x,F)=0\}.$$
\textup{(}The definition of $N(F)$ does not depend on $x$\/\textup{)}. Define
$$h_F(x):=P_x(N(F)).$$
Then $0\leq h_F(x)\leq 1$, $R_Wh_F=h_F$, $h_F$ is continuous and for $P_x$-a.e.\ $\omega\in\Omega$
$$\lim_{n\rightarrow\infty}h_F(\tau_{\omega_n}\cdots \tau_{\omega_1}x)=\left\{\begin{array}{cc}
1,&\mbox{ if }\omega\in N(F),\\
0,&\mbox{ if }\omega\not\in N(F).\end{array}\right.$$
\end{proposition}

\begin{proof}
Since the maps $\tau_l$ are contractions, it follows that $$\lim_nd(\tau_{\omega_n}\cdots \tau_{\omega_1}x,\tau_{\omega_n}\cdots \tau_{\omega_1}y)=0$$ for all $x,y\in\br^d$; hence the definition of $N(F)$ does not depend on $x$. 
\par
Consider the characteristic function $V_F(x,\omega):=\chi_{N(F)}(\omega)$, $x\in\br^d$, $\omega\in\Omega$. Then 
$$V_F(x,\omega_1\omega_2\wdots )=V_F(\tau_{\omega_1}x,\omega_2\omega_3\wdots ).$$
And $h_F(x)=P_x(V_F(x,\cdot))$.
The previous discussion in Section \ref{pathmeasures} then proves all the statements in the proposition.
\end{proof}
\par
In conclusion, this shows that every invariant set $F$ comes along with a naturally associated harmonic function $h_F$ ; see also Lemma \ref{lemhaco} above.

\begin{proposition}\label{propfinite}\cite[Propostion 2.3]{CCR96}
There exists a constant $\delta>0$ such that for any two disjoint compact invariant subsets $F$ and $G$, $d(F,G)>\delta$. There is only a finite number of minimal compact invariant subsets.
\end{proposition}
\begin{proof}
The first statement is in \cite{CCR96}. The only extra argument needed here is to prove that a minimal compact invariant subset is contained in some fixed compact set $K$. There is a norm which makes $S^{-1}$ a contraction. Define 
$K$ to be the closed ball centered at the origin with radius 
$$\rho:=\sup_{l\in L}\|l\|\frac{\|S^{-1}\|}{1-\|S^{-1}\|}.$$
Then $K$ is invariant for all maps $\tau_l$, $l\in L$, and 
$$\lim_{n\rightarrow\infty}d(\tau_{\omega_n}\cdots \tau_{\omega_1}x,K)=0\quad(x\in\br^d,\omega\in\Omega).$$
(See \cite[page 163]{CCR96}).
\par
If $F$ is a minimal compact invariant subset then take $x\in F$, and take $y$ to be one of the accumulation points of one of the trajectories. Then $y\in F\cap K$. With Proposition \ref{propmini}, $F\subset K$. The second statement follows.
\end{proof}

\begin{proposition}\label{propsupp}
Let $F_1, F_2,\dots , F_p$ be a family of mutually disjoint closed invariant subsets of $\br^d$ such that there is no closed invariant set $F$ with $F\cap\bigcup_kF_k=\emptyset$. Then 
$$P_x\left(\bigcup_{k=1}^pN(F_k)\right)=1\quad(x\in\br^d).$$
\end{proposition}

\begin{proof}
We reason by contradiction. Assume that for some $x\in\br^d$, $P_x(\bigcup_kN(F_k))<1$. Then define the function 
$$h(x)=P_x\left(\bigcup_kN(F_k)\right)=\sum_{k=1}^ph_{F_k}(x)<1.$$
According to Proposition \ref{propnf}, $R_Wh=h$ and $h$ is continous.
\par
Using again Proposition \ref{propnf}, there are some paths $\omega\not\in\bigcup_kN(F_k)$ such that 
$$\lim_{n\rightarrow\infty}h(\tau_{\omega_n}\cdots \tau_{\omega_1}x)=0.$$
Since $h$ is continous this implies that the set $Z$ of the zeroes of $h$ is not empty.
The equation $R_Wh=h$ also shows that $Z$ is a closed invariant subset.
\par
We show that $Z$ is disjoint from $\bigcup_k F_k$. If $Z\cap F_k\neq\emptyset$ for some $k\in\{1,\dots ,p\}$ then take $y\in F_k\cap Z$. There exists $\omega\in\Omega$ such that $W_B(\tau_{\omega_n}\cdots \tau_{\omega_1}y)\neq 0$ for all $n\geq 1$. (This is because $\sum_{l\in L}W_B(\tau_lz)=1$ for all $z$, so a transition is always possible.) But then, by invariance, $\tau_{\omega_n}\cdots \tau_{\omega_1}y\in F_k\cap Z$. This implies $\omega\in N(F_k)$ so, by Proposition \ref{propnf}, $\lim_nh_{F_k}(\tau_{\omega_n}\cdots \tau_{\omega_1}x)=1$. On the other hand $\tau_{\omega_n}\cdots \tau_{\omega_1}y\in Z$ so $h(\tau_{\omega_n}\cdots \tau_{\omega_1}y)=0$ for all $n\geq 1$. This yields the contradiction.
\par
Thus $Z$ is disjoint from $\bigcup_k F_k$, and this contradicts the hypothesis, and the proposition is proved.
\end{proof}

\begin{remark}
A family $F_1,\dots ,F_p$ as in Proposition \ref{propsupp} always exists because one can take all the minimal compact invariant sets. Proposition \ref{propfinite} shows that there are only finitely many such sets. And since every closed invariant set contains a minimal one, this family will satisfy the requirements. 
\end{remark}

\begin{theorem}\label{thcora}\cite[Th\'eor\`eme 2.8]{CCR96}
Let $M$ be minimal compact invariant set contained in the set of zeroes of an entire function $h$ on $\br^d$.
\begin{enumerate}
\item[a)] There exists $V$, a proper subspace of $\br^d$ invariant for $S$ \textup{(}possibly reduced to $\{0\}$\textup{)}, such that $M$ is contained in a finite union $\mathcal{R}$ of translates of $V$.
\item[b)] This union contains the translates of $V$ by the elements of a cycle\\  $\{x_0, \tau_{l_1}x_0,\dots , \tau_{l_{m-1}}\cdots \tau_{l_1}x_0 \}$ contained in $M$, and for all $x$ in this cycle, the function $h$ is zero on $x+V$.
\item[c)] Suppose the hypothesis ``(H) modulo $V$'' is satisfied, i.e., for all $p\geq0$ the equality $\tau_{k_1}\cdots \tau_{k_p}0-\tau_{k_1'}\cdots \tau_{k_p'}0\in V$, with $k_i,k_i'\in L$ implies $k_i-k_i'\in V$ for all $i\in\{1,\dots ,p\}$. Then 
$$\mathcal{R}=\{x_0+V,\tau_{l_1}x_0+V,\dots ,\tau_{l_{m-1}}\cdots \tau_{l_1}x_0+V\},$$
and every possible transition from a point in $M\cap\tau_{l_q}\cdots \tau_{l_1}x_0+V$ leads to a point in $M\cap\tau_{l_{q+1}}\cdots \tau_{l_1}x_0+V$ for all $1\leq q\leq m-1$, where $\tau_{l_m}\cdots \tau_{l_1}x_0=x_0$.
\item[d)] Since the function $W_B$ is entire, the union $\mathcal{R}$ is itself invariant.
\end{enumerate}
\end{theorem}
\par
A particular example of a minimal compact invariant set is a $W_B$-cycle. In this case, the subspace $V$ in Theorem \ref{thcora} can be take to be $V=\{0\}$:
\begin{definition} 
A {\it cycle} of length $p$ for the IFS $(\tau_l)_{l\in L}$ is a set of (distinct) points of the form $\mathcal{C}:=\{x_0,\tau_{l_1}x_0,\dots ,\tau_{l_{m-1}}\cdots \tau_{l_1}x_0\}$, such that $\tau_{l_m}\cdots \tau_{l_1}x_0=x_0$, with $l_1,\dots ,l_m\in L$. A $W_B$-{\it cycle} is a cycle $\mathcal{C}$ such that $W_B(x)=1$ for all $x\in\mathcal{C}$.
\par
For a finite sequence $l_1,\dots ,l_m\in L$ we will denote by $\underline{l_1\wdots l_m}$ the path in $\Omega$ obtained by an infinite repetition of this sequence
$$\underline{l_1\wdots l_m}:=(l_1\wdots l_ml_1\wdots l_m\wdots )$$
\end{definition}

\section{Statement of results}
\par
In the next definition we describe a way a given affine IFS $(\br^d, (\tau_b)_{b\in B})$, might factor such that the Hadamard property of Definition \ref{defhada} is preserved for the two factors. As a result we get a notion of reducibility (Definition \ref{defredcond}) for this class of affine IFSs. 
\begin{definition}\label{defreducible}
We say that the Hadamard triple $(R,B,L)$ is {\it reducible} to $\br^r$ if the following conditions are satisfied
\begin{enumerate}
\item\label{defreduciblei}
The subspace $\br^r\times\{0\}$ is invariant for $R^T$, so $S=R^T$ has the form 
$$S=\left[\begin{array}{cc}
S_1&C\\
0&S_2\end{array}\right]\quad
S^{-1}=\left[\begin{array}{cc}
S_1^{-1}&D\\
0&S_2^{-1}\end{array}\right],$$
with $S_1,C, S_2$ integer matrices.
\item\label{defreducibleii}
The set $B$ has the form $\{(r_i,\eta_{i,j})\,|\,i\in\{1,\dots ,N_1\},j\in\{1,\dots ,N_2\}\}$ where $r_i$ and $\eta_{i,j}$ are integer vectors;
\item\label{defreducibleiii}
The set $L$ has the form $\{(\gamma_{i,j},s_j)\,|\,j\in\{1,\dots ,N_2\},i\in\{1,\dots ,N_1\}\}$ where $s_j,\gamma_{i,j}$ are integer vectors;
\item\label{defreducibleiv}
$(S_1^T,\{r_i\,|\, i\in\{1,\dots ,N_1\}, \{\gamma_{i,j}\,|\,i\in\{1,\dots ,N_1\}\})$ is a Hadamard triple for all $j$;
\item\label{defreduciblev}
$(S_2^T,\{\eta_{i,j}\,|\,i\in\{1,\dots ,N_2\}\},\{s_j\,|\,j\in\{1,\dots ,N_2\}\})$ is a Hadamard triple for all $i$;
\item\label{defreduciblevi}
The invariant measure for the iterated function system $$\tau_{r_i}(x)=(S_1^T)^{-1}(x+r_i)\quad(x\in\br^r),i\in\{1,\dots ,N_1\}$$ is a spectral measure, and has no overlap, i.e., $\mu_1(\tau_{r_i}(X_1)\cap\tau_{r_j}(X_1))=0$ for all $i\neq j$, where $X_1$ is the attractor of the IFS $(\tau_{r_i})_{i\in\{1,\dots ,N_1\}}$.
\end{enumerate}
For convenience we will allow $r=0$, and every Hadamard triple is trivially reducible to $\br^0=\{0\}$. Note also that these conditions imply that $N=N_1N_2$.
\end{definition}

\begin{proposition}\label{propredu}
Let $(R,B,L)$ be a Hadamard triple such that $\br^r\times\{0\}$ is invariant for $R^T$. Assume that for all $b_1\in\proj_{\br^r}(B)$, the number of $b_2\in\br^{d-r}$ such that $(b_1,b_2)\in B$ is $N_2$, independent of $b_1$, and for all $l_2\in\proj_{\mathbb{R}^{d-r}}(L)$, the number of $l_1\in\br^r$ such that $(l_1,l_2)\in L$ is $N_1$, independent of $l_2$. Also assume that $N_1N_2=N$. Then the conditions \textup{(\ref{defreduciblei})--(\ref{defreduciblev})} in Definition \textup{\ref{defreducible}} are satisfied.
\end{proposition}

\begin{proof}
We define $\{r_1,\dots ,r_{M_1}\}:=\proj_{\br^r}(B)$. Using the assumption, for each $i\in\{1,\dots ,M_1\}$, we define 
$\{\eta_{i,1},\dots ,\eta_{i,N_2}\}$ to be the points in $\br^{d-r}$ with $(r_i,\eta_{i,j})\in B$. Similarly we can define $\{s_1,\dots ,s_{M_2}\}$, $\gamma_{i,j}$ for $L$. Since $M_1N_2=M_2N_1=N_1N_2=N$ we get $N_1=M_1$, $M_2=N_2$.
\par
Since the rows of the matrix $(e^{2\pi i R^{-1}b\cdot l})_{b\in B,l\in L}$ corresponding to $(r_{i_1},\eta_{i_1,j_1})$ and $(r_{i_1},\eta_{i_1,j_2})$ are orthogonal when $j_1\neq j_2$, and $i_1$ is fixed, we obtain (with the notation in Definition \ref{defreducible}):
$$\sum_{i=1}^{N_1}\sum_{j=1}^{N_2}e^{2\pi i (\eta_{i_1,j_1}-\eta_{i_1,j_2})\cdot S_2^{-1}s_j}=0,$$
and this implies (after dividing by $N_1$) that the rows of the matrix $(e^{2\pi i \eta_{i_1,j}\cdot S_2^{-1}s_{j'}})$, with $j,j'\in\{1,\dots ,N_2\}$. are orthogonal. This proves (\ref{defreduciblev}). The statement in (\ref{defreducibleiv}) is obtained using the dual argument 
(use the transpose of $R$ and interchange $B$ and $L$).
\end{proof}
\begin{definition}
We say that two Hadamard triples $(R_1,B_1,L_1)$ and $(R_2,B_2,L_2)$ are {\it conjugate} if there exists a matrix $M\in GL_d(\bz)$ (i.e., $M$ is invertible, and $M$ and $M^{-1}$ have integer entries) such that $R_2=MR_1M^{-1}$, $B_2=MB_1$ and $L_2=(M^T)^{-1}L_1$.
\end{definition}

If the two systems are conjugate then the transition between the IFSs $(\tau_b)_{b\in B_1}$ and $(\tau_{Mb})_{b\in B_1}$ is done by the matrix $M$; and the transition betweeen the IFSs $(\tau_l)_{l\in L_1}$ and $(\tau_{(M^T)^{-1}l})_{l\in L_1}$ is done by the matrix $(M^T)^{-1}$.
\begin{proposition}\label{prop2_8}
If $(R_1,B_1,L_1)$ and $(R_2,B_2,L_2)$ are conjugate through the matrix $M$, then
\begin{enumerate}
\item\label{prop2_8(1)}
$\tau_{Mb_1}(Mx)=M\tau_{b_1}(x)$, $\tau_{(M^T)^{-1}l_1}((M^T)^{-1}x)=(M^T)^{-1}\tau_{l_1}(x)$, for all $b_1\in B_1$, $l_1\in L_1$;
\item\label{prop2_8(2)}
$W_{B_2}(x)=W_{B_1}(M^Tx)$ for all $x\in\br^d$;

\item\label{prop2_8(3)}
For the Fourier transform of the corresponding invariant measures, the following relation holds:
$\hat\mu_{B_2}(x)=\hat\mu_{B_1}(M^Tx)$ for all $x\in\br^d$;

\item\label{prop2_8(4)}
The associated path measures satsify the following relation:
$$P_x^2(E)=P_{M^Tx}^1(\{(M^Tl_1,M^Tl_2,\dots )\,|\,(l_1,l_2,\dots )\in E\}).$$
\end{enumerate}
\end{proposition}

\begin{definition}\label{defreducing}
Let $(R,B,L)$ be a Hadamard triple. We call a subspace $V$ of $\br^d$ {\it reducing} if there exists a Hadamard triple $(R',B',L')$, conjugate to $(R,B,L)$, which is reducible to $\br^r$, and such that the conjugating matrix $M$, i.e., $R'=MRM^{-1}$, maps $V$ onto $\br^r\times\{0\}$. We allow here $V=\{0\}$, and the trivial space is clearly reducing.
\end{definition}

\begin{definition}\label{defredcond}
We say that the Hadamard triple $(R,B,L)$ satisfies {\it the reducibility condition} if for all minimal compact invariant subsets $M$, the subspace $V$ given in Theorem \ref{thcora} can be chosen to be reducing, and, for any two distinct  minimal compact invariant sets $M_1$, $M_2$, the corresponding unions $\mathcal{R}_1$, $\mathcal{R}_2$ of the translates of the associated subspaces given in Theorem \ref{thcora} are disjoint. 
\end{definition}

\begin{proposition}\label{prophmodv}
If $V$ is a reducing subspace then the hypothesis ``(H) modulo $V$'' is satisfied. 
\end{proposition}
\begin{proof}
By conjugation we can assume $V=\br^r\times\{0\}$. We use the notations in 
Definition \ref{defreducible}. 
\par
Let $k_n,k_n'\in L$, $n\in\{1,\dots ,p\}$ such that 
$\tau_{k_1}\cdots \tau_{k_p}0-\tau_{k_1'}\cdots \tau_{k_p'}0\in V$. Then we can write $k_n=(\eta_{i_n,j_n},s_{j_n})$,
$k_n'=(\eta_{i_n',j_n'},s_{j_n'})$ for all $n\in\{1,\dots ,p\}$. Then by a computation we obtain
$$\sum_{n=1}^pS_2^{-n}(s_{j_n}-s_{j_n'})=0.$$
This implies
$$\sum_{n=1}^pS_2^{p-n}(s_{j_n}-s_{j_n'})=0.$$
However, the Hadamard condition (\ref{defreduciblev}) in Definition \ref{defreducible} implies, according to Remark \ref{reminco}, that $s_{j_p}$ and $s_{j_p'}$ are not congruent $\mod S_2$, unless $j_p=j_p'$. Thus $j_p=j_p'$. By induction we obtain that $j_n=j_n'$ for all $n$ and this implies the hypothesis ``(H) modulo $V$''.
\end{proof}
\begin{theorem}\label{thmain}
Let $R$ be an expanding $d\times d$ integer matrix, $B$ a subset of $\bz^d$ with $0\in B$. Assume that there exists a subset $L$ of $\bz^d$ with $0\in L$ such that $(R,B,L)$ is a Hadamard triple which satisfies the reducibility condition. Then the invariant measure $\mu_B$ is a spectral measure.
\end{theorem}

\begin{remark}
If for all minimal compact invariant sets one can take the subspace $V$ to be $\{0\}$, i.e., if all the minimal compact invariant subsets are $W_B$-cycles, then the reducibility condition is automatically satisfied, and we reobtain Theorem 7.4 from \cite{DuJo05}.
\end{remark}

\section{Proofs}\label{sect4}
\par
The idea of the proof is to use the relation $\sum_{F}h_F=1$ from Proposition \ref{propsupp}. The functions $h_F$ will be written in terms of $|\hat\mu_B|^2$, and this relation will translate into the Parseval equality for a family of exponential function.

\addvspace{\medskipamount}\noindent\textbf{Invariant sets and invariant subspaces.}
We want to evaluate first $h_F(x)=P_x(N(F))$ for minimal invariant sets $F$. Theorem \ref{thcora} will give us the structure of these sets and this will aid in the computation.
\par
Consider a minimal compact invariant set $M$. Using Theorem \ref{thcora} we can find an invariant subspace $V$ such that $M$ is contained in the union of some translates of $V$. Since the reducibility condition is satisfied, we can take $V$ reducible. Proposition \ref{prophmodv} shows that the hypothesis ``(H) modulo $V$'' is satisfied. Therefore we can use part (c) of the theorem, and conclude that, for some cycle
$\mathcal{C}:=\{x_0,\tau_{l_1}x_0,\dots ,\tau_{l_{m-1}}\cdots \tau_{l_1}x_0\}$, with $\tau_{l_m}\cdots \tau_{l_1}x_0=x_0$, $M$ is contained in the union 
$$\mathcal{R}=\{x_0+V,\tau_{l_1}x_0+V,\dots ,\tau_{l_{m-1}}\cdots \tau_{l_1}x_0+V\},$$
and $\mathcal{R}$ is an invariant subset. 
\par
By conjugation we can assume first that $V=\br^r\times\{0\}$, and the Hadamard triple $(R,B,L)$ is reducible to $\br^r$.
We will use the notations in Definition \ref{defreducible}. Thus $S$, $B$ and $L$ have the specific form given in this definition. Also, points in $\br^d$ are of the form $(x,y)$ with $x\in\br^r$ and $y\in\br^{d-r}$. We refer to $x$ as the ``first component'' and to $y$ as the ``second component''.
For a path $(\omega_1\wdots \omega_k\wdots )$ in $\Omega$ we will use the notation $(\omega_{1,1}\wdots \omega_{k,1}\wdots )$ for the path of the first components, and $(\omega_{1,2}\wdots \omega_{k,2}\wdots )$ for the path of the second components.
\par
We will also consider the IFS defined on the second component:
$$\tau_{s_i}(y)=S_2^{-1}(y+s_i)\quad(y\in\br^{d-r},i\in\{1,\dots ,N_2\}).$$

\par
We want to compute $P_{(x,y)}(N(\mathcal{R}))$ (see Proposition \ref{propnf} for the definition of $N(\mathcal{R})$).
\begin{lemma}\label{lem2_1}
Let $h_1,\dots ,h_m\in\{s_i\,|\,i\in\{1,\dots ,N_2\}\}$ be the second components of the sequence $l_1,\dots ,l_m$ that defines the cycle $\mathcal{C}$. A path $(\omega_1\omega_2\wdots )$ is in $N(\mathcal{R})$ if and only if the second component of this path is of the form $(\omega_{1,2}\wdots \omega_{k,2}\underline{h_1\wdots h_m})$, where $\omega_{1,2},\dots ,\omega_{k,2}$ are arbitrary in $\{s_i\,|\,i\in\{1,\dots ,N_2\}\}$.
\end{lemma}

\begin{proof}
Since $V=\br^r\times\{0\}$, the path $\omega$ is in $N(\mathcal{R})$ if and only if the second component of $\tau_{\omega_n}\cdots \tau_{\omega_1}(x,y)$ approaches the set $\mathcal{C}_2$ of the second components of the cycle $\mathcal{C}$.  But note that 
$\tau_{(\omega_{k,1},\omega_{k,2})}(x,y)$ has the form $(*,\tau_{\omega_{k,2}}y)$. Thus we must have 
\begin{equation}\label{eqc2}
\lim_{n}d(\tau_{\omega_{k,2}}\cdots \tau_{\omega_{1,2}}y,\mathcal{C}_2)= 0.
\end{equation}
Also $\mathcal{C}_2=\{y_0,\tau_{h_1}y_0,\dots ,\tau_{h_{m-1}}\cdots \tau_{h_1}y_0\}$ is a cycle for the IFS $(\tau_{s_i})_i$, where $y_0$ is the second component of $x_0$, and $\tau_{h_m}\cdots \tau_{h_1}y_0=y_0$.  But then (\ref{eqc2}) is equivalent to the fact that the path $(\omega_{1,2}\omega_{2,2}\wdots )$ ends in an infinite repetition of the cycle $h_1\wdots h_m$ (see \cite[Remark 6.9]{DuJo05}). This proves the lemma.
\end{proof}
\par
Thus the paths in $N(\mathcal{R})$ are arbitrary on the first component, and end in a repetition of the cycle on the second. We will need to evaluate the following quantity, for a fixed $l_2\in\{s_1,\dots ,s_{N_2}\}$, and $(x,y)\in\br^d$:
\begin{align*}
A&:=\sum_{l_1\mbox{ with }(l_1,l_2)\in L}W_B(\tau_{(l_1,l_2)}(x,y))
\\
&\phantom{:}=\sum_{l_1}\frac{1}{N_1^2N_2^2}\sum_{i,i'}\sum_{j,j'}e^{2\pi i((r_i-r_{i'})\cdot(S_1^{-1}(x+l_1)+D(y+l_2))+(\eta_{i,j}-\eta_{i',j'})\cdot(S_2^{-1}(y+l_2)))}.
\end{align*}
But, because of the Hadamard property (\ref{defreducibleiv}) in Definition \ref{defreducible}, 
$$\frac{1}{N_1}\sum_{l_1}e^{2\pi i(r_i-r_{i'})\cdot S_1^{-1}l_1}=\left\{\begin{array}{cc}
1,&i=i',\\
0,&i\neq i'.\end{array}\right.$$
Therefore
$$A=\frac{1}{N_1N_2^2}\sum_i\sum_{j,j'} e^{2\pi i(\eta_{i,j}-\eta_{i,j'})\cdot S_2^{-1}(y+l_2)}$$
and
\begin{equation}\label{eqwtilde}
\sum_{l_1}W_B(\tau_{(l_1,l_2)}(x,y))=\frac{1}{N_1}\sum_{i=1}^{N_1}W_i(\tau_{l_2}y)=:\tilde W(\tau_{l_2}y),
\end{equation}
where
\begin{equation}\label{eqWi}
W_i(y)=\left|\frac{1}{N_2}\sum_{j=1}^{N_2}e^{2\pi i\eta_{i,j}\cdot y}\right|^2
\end{equation}

Next we compute $P_{(x,y)}$ for those paths that have a fixed second component $(l_{1,2}l_{2,2}\wdots l_{n,2}\wdots )$.
\begin{lemma}\label{lempatom}
$$P_{(x,y)}(\{(\omega_1\wdots \omega_n\wdots )\,|\,\omega_{n,2}=l_{n,2}\mbox{ for all }n\})=\prod_{k=1}^\infty\tilde W(\tau_{l_{k,2}}\cdots \tau_{l_{1,2}}y).$$
\end{lemma}

\begin{proof}
We compute for all $n$, by summing over all the possibilities for the first component, and using (\ref{eqcyl}):
\begin{multline*}
P_{(x,y)}(\{(\omega_1\omega_2\wdots )\,|\,\omega_{k,2}=l_{k,2}, 1\leq k\leq n\})
\\
=
\sum_{l_{1,1},\dots ,l_{n,1}}\prod_{k=1}^nW_B(\tau_{(l_{k,1},l_{k,2})}\cdots \tau_{(l_{1,1},l_{1,2})}(x,y))
=(*).
\end{multline*}
Using (\ref{eqwtilde}) we obtain further
\begin{multline*}
(*)=\tilde W(\tau_{l_{n,2}}\cdots \tau_{l_{1,2}}y)\sum_{l_{1,1},\dots ,l_{n-1,1}}\prod_{k=1}^{n-1}W_B(\tau_{(l_{k,1},l_{k,2})}\cdots \tau_{(l_{1,1},l_{1,2})}(x,y))
\\
=\dots 
=\prod_{k=1}^n\tilde W(\tau_{l_{k,2}}\cdots \tau_{l_{1,2}}y).
\end{multline*}
Then, letting $n\rightarrow\infty$ we obtain the lemma.
\end{proof}

Next we will see how the invariant measure $\mu_B$ and the attractor $X_B$ can be decomposed through the invariant subspace $V=\br^r\times\{0\}$.
\par
The matrix $R$ has the form:
\[
R=\left[\begin{array}{cc}
A_1&0\\
C^*&A_2\end{array}\right],\mbox{ and }
R^{-1}=\left[\begin{array}{cc}
A_1^{-1}&0\\
-A_2^{-1}C^*A_1^{-1}&A_2^{-1}\end{array}\right].
\]
By induction,
\[
R^{-k}=\left[\begin{array}{cc}
A_1^{-k}&0\\
D_k&A_2^{-k}\end{array}\right],\mbox{ where }D_k:=-\sum_{l=0}^{k-1}A_2^{-(l+1)}C^*A_1^{-(k-l)}.\]
We have 
\[
X_B=\{\sum_{k=1}^\infty R^{-k}b_k\,|\,b_k\in B\}.
\]
Therefore any element $(x,y)$ in $X_B$ can be written in the following form:
\[
x=\sum_{k=1}^\infty A_1^{-k}r_{i_k},\quad y=\sum_{k=1}^\infty D_kr_{i_k}+\sum_{k=1}^\infty A_2^{-k}\eta_{i_k,j_k}.
\]
\par

Define 
\[
X_1:=\{\sum_{k=1}^\infty A_1^{-k}r_{i_k}\,|\,i_k\in\{1,\dots ,N_1\}\}.
\]
Let $\mu_1$ be the invariant measure for the iterated function system 
\[
\tau_{r_i}(x)=A_1^{-1}(x+r_i),\quad i\in\{1,\dots ,N_1\}.
\]
The set $X_1$ is the attractor of this iterated function system.
\par
For each sequence $\omega=(i_1i_2\wdots )\in\{1,\dots ,N_1\}^{\bn}$, define $x(\omega)=\sum_{k=1}^\infty A_1^{-k}r_{i_k}$. Also, because of the non-overlap condition, for $\mu_1$-a.e.\ $x\in X_1$, there is a unique $\omega$ such that 
$x(\omega)=x$. We define this as $\omega(x)$. This establishes an a.e.\ bijective correspondence between $\Omega_1$ and $X_1$, $\omega\leftrightarrow x(\omega)$.
\par
Denote by $\Omega_1$ the set of all paths $(i_1i_2\wdots i_n\wdots )$ with $i_k\in\{1,\dots ,N_1\}$. For $\omega=(i_1i_2\wdots )\in\Omega_1$ define
\[
\Omega_2(\omega):=\{\eta_{i_1,j_1}\eta_{i_2,j_2}\wdots \eta_{i_n,j_n}\wdots \,|\,j_k\in\{1,\dots ,N_2\}\}.
\]
\par
For $\omega\in\Omega_1$ define $g(\omega):=\sum_{k=1}^\infty D_kr_{i_k}$, and $g(x):=g(\omega(x))$. Also we denote 
$\Omega_2(x):=\Omega_2(\omega(x))$.
\par
For $x\in X_1$, define
\[
X_2(x):=X_2(\omega(x)):=
\left\{\sum_{k=1}^\infty A_2^{-k}\eta_{i_k,j_k}\biggm|j_k\in\{1,\dots ,N_2\}\mbox{ for all }k\right\}.
\]
\par

Note that the attractor $X_B$ has the following form:
\[
X_B=\{(x,g(x)+y)\,|\,x\in X_1,y\in X_2(x)\}.
\]
We will show that the measure $\mu_B$ can also be decomposed as a product between the measure $\mu_1$ and some measures $\mu_{\omega}^2$ on $X_2(\omega)$.
\par
On $\Omega_2(\omega)$, consider the product probability measure $\mu(\omega)$ which assigns to each $\eta_{i_k,j_k}$ equal probabilities $1/N_2$. 
\par
Next we define the measure $\mu_{\omega}^2$ on $X_2(\omega)$. Let 
$r_\omega:\Omega_2(\omega)\rightarrow X_\omega^2$,
\[
r_\omega(\eta_{i_1,j_1}\eta_{i_2,j_2}\wdots )=\sum_{k=1}^\infty A_2^{-k}\eta_{i_k,j_k}.
\]
\par

Define the measure $\mu_x^2:=\mu_{\omega(x)}^2:=\mu_{\omega(x)}\circ r_{\omega(x)}^{-1}$.

\begin{lemma}\label{lem4_3}
Let $\sigma$ be the shift on $\Omega_1$, $\sigma(i_1i_2\wdots )=(i_2i_3\wdots )$. Let $\omega=(i_1i_2\wdots )\in\Omega_1$. 
Then for all measurable sets $E$ in $X_2(\omega)$,
\[
\mu_\omega^2(E)=\frac{1}{N_2}\sum_{j=1}^{N_2}\mu_{\sigma(\omega)}^2(\tau_{\eta_{i_1,j}}^{-1}(E)).
\]
The Fourier transform of the measure $\mu_\omega^2$ satisfies the equation:
\begin{equation}\label{eqmu2}
\hat\mu_\omega^2(y)=m(S_2^{-1}y,i_1)\hat\mu_{\sigma(\omega)}^2(S_2^{-1}y),
\end{equation}
where
$$m(y,i_1)=\frac{1}{N_2}\sum_{j=1}^{N_2}e^{2\pi i\eta_{i_1,j}\cdot y}.$$

\end{lemma}

\begin{proof}
We define the maps $\xi_{\eta_{i_1,j}}:\Omega_2(\sigma(\omega))\rightarrow\Omega_2(\omega)$, $$\xi_{\eta_{i_1,j}}(\eta_{i_2,j_2}\eta_{i_3,j_3}\wdots )=
(\eta_{i_1,j_1}\eta_{i_2,j_2}\wdots ).$$
\par
Then $r_\omega\circ\xi_{\eta_{i_1,j}}=\tau_{\eta_{i_1,j}}\circ r_{\sigma(\omega)}$.
\par
The relation given in the lemma can be pulled back through $r_\omega$ to the path spaces $\Omega_2(\omega)$, and becomes equivalent to:
$$\mu_\omega(E)=\frac{1}{N_2}\sum_{j}\mu_{\sigma(\omega)}(\xi_{\eta_{i_1,j}}^{-1}(E)),$$ 
and this can be immediately be verified on cylinder sets, i.e., the sets of paths in $\Omega_2(\omega)$ with some prescribed first $n$ components. 
\par
{}From this it follows that 
\[
\int f\,d\mu_\omega^2=\frac{1}{N_2}\sum_{j=1}^{N_2}\int f\circ\tau_{\eta_{i_1,j}}\,d\mu_{\sigma(\omega)}^2.
\]
Applying this to the function $s\mapsto e^{2\pi i s\cdot y}$ we obtain equation (\ref{eqmu2}).
\end{proof}

\begin{lemma}\label{lem4_4}
\[
\int_{X_B}f\,d\mu_B=\int_{X_1}\int_{X_2(x)}f(x,y+g(x))\,d\mu_x^2(y)\,d\mu_1(x).
\]
\end{lemma}
\begin{proof}
We begin with a relation for the function $g$.
\begin{equation}
g(A_1^{-1}(x+r_i))=D_1(x+r_i)+A_2^{-1}g(x)
\end{equation}

Indeed, if $\omega(x)=(i_1i_2\wdots )$, then $\omega(A_1^{-1}(x+r_i))=(ii_1i_2\wdots )$.
So
\begin{align*}
g(A_1^{-1}(x+r_i))&=D_1r_i+\sum_{k=1}^\infty D_{k+1}r_{i_k}=
D_1r_i-\sum_{k=1}^\infty\sum_{l=0}^{k}A_2^{-(l+1)}C^*A_1^{-(k+1-l)}r_{i_k}
\\
&=D_1r_i-\sum_{k=1}^\infty A_2^{-1}C^*A_1^{-k-1}r_{i_k}-\sum_{k=1}^\infty\sum_{l=0}^{k-1} A_2^{-(l+2)}C^*A_1^{-(k-l)}r_{i_k}
\\
&=D_1r_i+D_1x+A_2^{-1}g(x).
\end{align*}
\par

Next we show that the measure $\mu_B$ has the given decomposition. We check the invariance of the decomposition. We denote by $i_1(x)$, the first component of $\omega(x)$, and $\sigma(x)$ is the point in $X_1$ that corresponds to $\sigma(\omega(x))$.
\begin{align*}
&\int_{X_1}\int_{X_2(x)}f(x,y+g(x))\,d\mu_x^2(y)\,d\mu_1(x)
\\
&\qquad=\frac{1}{N_2}\sum_{j=1}^{N_2}\int_{X_1}\int_{X_2(\sigma(x))}f(x,A_2^{-1}(y+\eta_{i_1(x),j})+g(x))\,d\mu_{\sigma(x)}(y)\,d\mu_1(x)
\\
&\qquad=\frac{1}{N_1N_2}\sum_{i=1}^{N_1}\sum_{j=1}^{N_2}\int_{X_1}\int_{X_2(\sigma(\tau_ix))}f(A_1^{-1}(x+r_i),A_2^{-1}(y+\eta_{i_1(\tau_{r_i}x),j})
\\
&\qquad\qquad\quad{}+g(A_1^{-1}(x+r_i)))\,d\mu_{\sigma(\tau_{r_i}x)}(y)\,d\mu_1(x)
\\
&\qquad=\frac{1}{N}\sum_{i,j}\int_{X_1}\int_{X_2(x)}f(A_1^{-1}(x+r_i),
\\
&\qquad\qquad\quad{}D_1(x+r_i)+A_2^{-1}(y+g(x)+\eta_{i,j}))\,d\mu_x^2(y)\,d\mu_1(x)
\\
&\qquad=\frac{1}{N}\sum_{i,j}\int_{X_1}\int_{X_2(x)}f\circ\tau_{(r_i,\eta_{i,j})}(x,y+g(x))\,d\mu_x^2(y)\,d\mu_1(x).
\end{align*}
Using the uniqueness of the invariant measure for an IFS, we obtain the lemma.
\end{proof}

\begin{lemma}\label{lem4_5}
If $\Lambda_1$ is a spectrum for the measure $\mu_1$, then 
$$F(y):=\sum_{\lambda_1\in\Lambda_1}|\hat\mu_B(x+\lambda_1,y)|^2=\int_{X_1}|\hat\mu_s^2(y)|^2\,d\mu_1(s)\quad(x\in\br^r,y\in\br^{d-r}).$$
\end{lemma}

\begin{proof}
\begin{align*}
F(y)&=\sum_{\lambda_1}\left|\int_{X_1}\int_{X_2(s)}e^{2\pi i((x+\lambda_1)\cdot s+y\cdot(t+g(s))}\,d\mu_x^2(t)\,d\mu_1(s)\right|^2
\\
&=\sum_{\lambda_1}\int_{X_1}\left(e^{2\pi i(x\cdot s+y\cdot g(s))}\hat\mu_x^2(y)\right)e^{2\pi i\lambda_1\cdot s}\,d\mu_1(s)
\\
&=\int_{X_1}|\hat\mu_s^2(y)|^2\,d\mu_1(s),
\end{align*}
where we used the Parseval identity in the last equality.

\end{proof}

\begin{lemma}
$$F(y)=\tilde W(S_2^{-1}y)F(S_2^{-1}y).$$
Also 
$$F(y)=\prod_{k=1}^\infty W(S_2^{-k}y)\quad(y\in\br^{d-r}).$$
\end{lemma}

\begin{proof}
Using Lemma \ref{lem4_3} and Lemma \ref{lem4_5}, and the fact that
$$\frac{1}{N_1}\sum_{i=1}^{N_1}|m(y,i)|^2=\tilde W(y)\quad(y\in\br^{d-r}),$$
we obtain
\begin{align*}
F(y)&=\int_{X_1}\left|m(S_2^{-1}y,i_1(s))\right|^2\left|\hat\mu_{\sigma(s)}^2(S_2^{-1}y)\right|^2\,d\mu_1(s)
\\
&=\frac{1}{N_1}\sum_{i=1}^{N_1}\int_{X_1}\left|m(S_2^{-1}y,i_1(\tau_{r_i}s)
)\right|^2\left|\hat\mu_{\sigma(\tau_{r_i}s)}^2(S_2^{-1}y)\right|^2\,d\mu_1(s)
\\
&=\vphantom{\int_{X_1}}\tilde W(S_2^{-1}y)F(S_2^{-1}y).
\end{align*}
We also have $F(0)=1$ because $\mu_1$ and $\mu_{\omega}^2$ are probability measures. Using Lemma \ref{lem4_5} it is easy to see that $F$ is continuous. Also $\tilde W(0)=1$ and for some $0<c<1$, $\|S_2^{-k}\|\leq c^k$ for all $k$ (because $S_2$ is expansive), and $\tilde W$ is Lipschitz, the infinite product is then convergent to $F(y)$.
\end{proof}

\par
Now consider the cycle associated to the minimal invariant set $M$,  $$\mathcal{C}=\{x_0,\tau_{l_1}x_0,\dots ,\tau_{l_{m-1}}\cdots \tau_{l_1}x_0\}$$ as described in the begining of the section, with $\tau_{l_m}\cdots \tau_{l_0}x_0=x_0$. Consider the second components of this cycle. Let the second component of $x_0$ be $y_0$ and let $h_1,\dots ,h_m\in\{s_i\,|\,i\in\{1,\dots ,N_2\}\}$ be the second components of $l_1,\dots ,l_m$. 
\begin{lemma}\label{lemwtildecycle}
The set $\mathcal{C}_2:=\{y_0,\tau_{h_1}y_0,\dots ,\tau_{h_{m-1}}\cdots \tau_{h_1}y_0\}$ is a $\tilde W$-cycle.
\end{lemma}

\begin{proof}
We saw in the proof of Lemma \ref{lem2_1} that $\mathcal{C}_2$ is a cycle. We only need to check that $\tilde W(y)=1$ for all $y\in\mathcal{C}_2$. Take the point $y_0$ and take some $s_j\neq h_1$. We claim that $\tau_{s_j}y_0$ cannot be one of the points in $\mathcal{C}_2$. Otherwise it would follow that $y_0$ is a fixed point for $\tau_{\omega_q}\cdots \tau_{\omega_1}$, for some $\omega_1,\omega_2,\dots ,\omega_q\in\{s_j\,|\,j\in\{1,\dots ,N_2\}\}$  with $\omega_1=s_j\neq h_1$. But $x_0$ is also a fixed point for $\tau_{h_m}\cdots \tau_{h_1}$. It follows that $x_0$ is fixed also by $(\tau_{h_m}\cdots \tau_{h_1})^q$ and $(\tau_{\omega_q}\cdots \tau_{\omega_1})^m$. Writing the corresponding fixed point equations, we obtain:
$$(S_2^{mq}-I)^{-1}(h_1+Sk)=x_0=(S_2^{mq}-I)^{-1}(\omega_1+Sk'),$$
for some $k,k'\in\bz^{d-r}$. But this implies that $h_1\equiv\omega_1\mod S_2\bz^{d-r}$ and this is impossible because of the Hadamard property (v) in Definition \ref{defreducible} and Remark \ref{reminco}. This proves our claim.
\par
Since $\tau_{s_j}y_0$ is not in $\mathcal{C}_2$, the invariance of the set $\mathcal{R}=\bigcup_{y\in\mathcal{C}_2} (y+\br^r\times\{0\})$ implies that $W_B(\tau_{(\eta_{i,j},s_j)}(x,y_0))=0$ for all $i\in\{1,\dots ,N_2\}$. 
But then, with equation (\ref{eqwtilde}), this implies that $\tilde W(\tau_{s_j}y_0)=0$, for all $s_j\neq h_1$. And since 
$$\sum_{j=1}^{N_2}\tilde W(\tau_{s_j}y_0)=1,$$
it follows that $\tilde W(\tau_{h_1}y_0)=1$. The same argument works for the other points in $\mathcal{C}_2$, and we obtain the result.
\end{proof}

\begin{lemma}\label{lemplxc}
The following relation holds for all $k\geq0$:
$$\tilde W(y+S_2^{km}y_0)=\tilde W(y)\quad(y\in\br^d).$$
\end{lemma}

\begin{proof}
Since $\tilde W(y_0)=1$, it follows that $W_i(y_0)=1$ for all $i\in\{1,\dots ,q_1\}$. Therefore all the terms in the sum which defines $W_i$ must be $1$ which means that $\eta_{i,j}\cdot y_0\in\bz$ for all $i,j$. This implies that 
$W_i(y+y_0)=W_i(y)$
\par
On the other hand, as $y_0$ is a fixed point for $\tau_{h_{m}}\cdots \tau_{h_1}$, we have $S_2^my_0
\equiv y_0\mod \bz^{d-r}.$ By induction $S_2^{km}y_0\equiv y_0\mod\bz^{d-r}$ for all $k\geq0$.
\par
Thus, $\tilde W(y+S_2^{km}y_0)=\tilde W(y+y_0)=\tilde W(y)$.
\end{proof}

\begin{lemma}\label{lem4_9}
For $\omega=\omega_0\wdots \omega_{km-1}\in\{s_j\,|\,j\in\{1,\dots ,N_2\}\}^{km}$, define
$E_{\omega,\mathcal{C}}$ to be the set of paths in $\Omega$ that have the second component equal to $(\omega_0\wdots \omega_{kp-1}\underline{h_1\wdots h_m})$, and
$$k_{\mathcal{C}}(\omega):=\omega_0+\dots +S_2^{km-1}\omega_{km-1}-S_2^{km}y_0.$$
Then
$$P_{(x,y)}(E_{\omega,C})=F(y+k_C(\omega))\quad(x\in\br^r,y\in\br^{d-r}).$$
\end{lemma}

\begin{proof}
For $q\leq km-1$
$$\tau_{\omega_{q-1}}\cdots \tau_{\omega_0}y=S_2^{-q}(y+\omega_0+\dots +S_2^{q-1}\omega_{q-1})$$
$$S_2^{-q}(y+k_C(\omega))\equiv S_2^{-q}(y+\omega_0+\dots +S_2^{q-1}\omega_{m-1})-S_2^{-q+km}y_0\mod\bz^d$$
But $S_2^{-q+km}y_0\equiv y'\mod\bz^d$ for one of the elements $y'$ of the $\tilde W$-cycle $\mathcal{C}_2$. 
\par
Therefore $\tilde W(\tau_{\omega_{q-1}}\cdots \tau_{\omega_0}y)=\tilde W(S_2^{-q}(y+k_C(\omega))$.
\par
Next, for $j\geq k$, 
\begin{align*}
&(\tau_{h_m}\cdots \tau_{h_1})^{j-k}\tau_{\omega_{km-1}}\cdots \tau_{\omega_0}y
\\
&\qquad=S_2^{-jm}(y+\omega_0+\dots +S_2^{km-1}\omega_{km-1}
\\
&\qquad\qquad\quad{}+S_2^{km}(I+S_2^m+\dots +S_2^{(j-k-1)m})(h_1+\dots +S_2^{m-1}h_{m}))=(*).
\end{align*}
Using $y_0=(S_2^{m}-I)^{-1}(h_1+\dots +S_2^{m-1}h_m)$,
\begin{align*}
(*)&=S_2^{-jm}(y+\omega_0+\dots +S_2^{km-1}\omega_{km-1}
\\
&\qquad\qquad\quad{}+S_2^{km}(S_2^{(j-k)m}-I)(S_2^m-I)^{-1}(S_2^m-I)y_0)
\\
&\qquad=S_2^{-jm}(y+k_C(\omega)+S_2^{jm}y_0).
\end{align*}
Using Lemma \ref{lemplxc}, we obtain that
$$\tilde W((\tau_{h_{m}}\cdots \tau_{h_1})^{j-k}\tau_{\omega_{km-1}}\cdots \tau_{\omega_0}y)=\tilde W(S_2^{-jm}(y+k_C(\omega))).$$
Also, using the previous equalities, for $q\leq m$,
$$\tau_{h_q}\cdots \tau_{h_1}(\tau_{h_m}\cdots \tau_{h_1})^{j-k}\tau_{\omega_{km-1}}\cdots \tau_{\omega_0}y=$$
$$=\tau_{h_q}\cdots \tau_{h_1}(S_2^{-jm}(y+k_{\mathcal{C}}(\omega))+y_0)=\tau_{h_q}\cdots \tau_{h_1}y_0+S_2^{-jm-q}(y+k_{\mathcal{C}}(\omega))$$
and since $\tau_{h_q}\cdots \tau_{h_1}$ is also an element of the $\tilde W$-cycle, Lemma \ref{lemplxc} applies and
$$\tilde W(\tau_{h_q}\cdots \tau_{h_1}(\tau_{h_{m}}\cdots \tau_{h_1})^{j-k}\tau_{\omega_{km-1}}\cdots \tau_{\omega_0}y)=\tilde W(S_2^{-jm-q}(y+k_{\mathcal{C}}(\omega))).$$
\par
This proves, using the infinite product formulas for $P_{(x,y)}$ and $F$ in Lemma \ref{lempatom} and Lemma \ref{lem4_5} that
$$P_{(x,y)}(E_{\omega,C})=F(y+k_{\mathcal{C}}(\omega)).$$

\end{proof}

\begin{proposition}\label{propnr}
There exists a set $\Lambda(M)\subset\bz^d$ such that 
$$h_{\mathcal{R}}(x)=P_x(N(\mathcal{R}))=\sum_{\lambda\in\Lambda(M)}|\hat\mu_B(x+\lambda)|^2\quad(x\in\br^d).$$
\end{proposition}

\begin{proof}
First note that, with proposition \ref{prop2_8}, we can assume that the Hadamard triple $(R,B,L)$ is reducible to $\br^r$ and 
$V=\br^r\times 0$. 
\par
With Lemma \ref{lem2_1} we see that $N(\mathcal{R})$ is the set of all paths such that the second component has the form $(\omega_0\wdots \omega_{k}\underline{s_{j_{0}}\wdots s_{j_{p-1}}})$.
\par
We have
$$P_{(x,y)}(N(\mathcal{R}))=\sum_{\omega}P_{(x,y)}(E_{\mathcal{C},\omega})$$
where the sum is indexed over all possible paths that end in a repetition of the cycle $h_1\wdots h_m$, so it can be indexed by a choice of a finite path $\omega_1\wdots \omega_{km-1}$ in with $\omega_i\in\{s_j\,|\,j\in\{1,\dots ,N_2\}\}^{km}$ for all $i$.
\par
Using Lemma \ref{lem4_9} and Lemma \ref{lem4_5} we obtain further: 
 $$P_{(x,y)}(N(\mathcal{R}))=\sum_{\omega}F(y+k_{\mathcal{C}}(\omega))$$
$$=\sum_{\omega}\sum_{\lambda\in\Lambda_1}|\hat\mu_B(x+\lambda_1,y+k_{\mathcal{C}}(\omega))|^2.$$
The proposition is proved. 
\end{proof}
\begin{remark}
It might happen that for two different paths $\omega$ the integers $k_{\mathcal{C}}$ are the same. Therefore the same $\lambda$ might appear twice in the set $\Lambda(M)$. We make the convention to count it twice. We will show in the end that actually this will not be the case.
\end{remark}
We are now in position to give the proof of the theorem. 

\begin{proof} (of Theorem \ref{thmain})
Let $M_1,\dots ,M_p$ be the list of all minimal compact invariant sets. The hypothesis shows that for each $k$ there is a reducing subspace $V_k$ and some cycle $\mathcal{C}_k$ such that $M_k\subset\mathcal{R}_k:=\mathcal{C}_k+V_k$, and moreover the sets $\mathcal{R}_k$ are mutually disjoint. With Proposition \ref{propnr} we see that there is some set $\Lambda(M_k)\subset\bz^d$ such that 
$$h_{\mathcal{R}_k}(x)=\sum_{\lambda\in\Lambda(M_k)}|\hat\mu_B(x+\lambda)|^2\quad(x\in\br^d).$$
With Proposition \ref{propsupp}, we have
\begin{equation}\label{eqparse}
1=\sum_{k=1}^ph_{\mathcal{R}_k}(x)=\sum_{k=1}^p\sum_{\lambda\in\Lambda(M_k)}|\hat\mu_B(x+\lambda)|^2.
\end{equation}
We check that a $\lambda$ cannot appear twice in the union of the sets $\Lambda(M_k)$. For some fixed $\lambda_0\in \bigcup_k\Lambda(M_k)$, take $x=-\lambda_0$ in (\ref{eqparse}). Since $\hat\mu_B(0)=1$, it follows that one term in the sum is $1$ (the one corresponding to $\lambda_0$) and the rest are $0$. Thus $\lambda_0$ cannot appear twice. 
Also for $\lambda\neq\lambda_0$, this implies that
$\hat\mu_B(-\lambda_0+\lambda)=0$ so the functions $e^{2\pi i\lambda_0\cdot x}$ and $e^{2\pi i\lambda\cdot x}$ are orthogonal in $L^2(\mu_B)$. 
\par
With the notation $e_x(t)=e^{2\pi i x\cdot t}$, we can rewrite (\ref{eqparse}) as 
$$\|e_{-x}\|_2^2=\sum_{\lambda\in\bigcup_{k=1}^p\Lambda(M_k)}|\ip{e_{-x}}{e_{\lambda}}|^2\quad(x\in\br^d).$$
But this, and the orthogonality, implies that the closed span of family of functions $\{e_{\lambda}\,|\,\lambda\in\Lambda\}$, where 
$\Lambda=\bigcup_{k=1}^p\Lambda(M_k)$, contains all functions $e_{x}$, and, by Stone-Weierstrass, this implies that it contains $L^2(\mu_B)$. Thus, $\{e_{\lambda}\,|\,\lambda\in\Lambda\}$ forms an orthonormal basis for $L^2(\mu_B)$.

\end{proof}

\section{Examples}\label{sect5}
Before we give the examples we will prove a lemma which helps in identifying candidates for the invariant subspaces containing minimal invariant sets. 
\begin{lemma}\label{lemposinv}
With the assumptions of Theorem \ref{thcora}, suppose that there is no proper subspace $W$ such that $X_B$ is contained in a finite union of translates of $W$. Let $V$ be an invariant subspace as in \ref{thcora}. Then there is some $x\in\br^d$ such that $W_B(x+v)=0$ for all $v\in V$. If in addition the hypothesis ``(H) modulo $V$'' is satisfied, and $\mathcal{C}:=\{x_0,\tau_{l_1}x_0,\dots ,\tau_{l_{m-1}}\cdots \tau_{l_1}x_0\}$ is the cycle given in Theorem \ref{thcora}, then $x$ can be taken to be any point $\tau_{l_k}\cdots \tau_{l_1}x_0$ of the cycle and $l$ can be taken to be any element of $L$ such that $l-l_{k+1}\not\in V$.
\end{lemma}

\begin{proof}
Consider the invariant union $\mathcal{R}$ of translates of $V$, as in Theorem \ref{thcora}. Then $\mathcal{R}$ cannot contain $X_B$ so for some $x\in\mathcal{R}$ and some $l\in L$ we have $\tau_l(x)\not\in\mathcal{R}$. But then, for all $v\in V$, $\tau_l(x+v)=\tau_lx+S^{-1}v$ cannot be in $\mathcal{R}$ (otherwise $\tau_lx=\tau_l(x+v)-S^{-1}v\in\mathcal{R}+V=\mathcal{R}$). Since $\mathcal{R}$ is invariant, it follows that $W_B(\tau_l(x+v))=0$. But $\tau_l(x+V)=\tau_lx+S^{-1}V=\tau_lx+V$ and this proves the first assertion.
\par
If $V$ also satisfies the hypothesis ``(H) modulo $V$'', then $\mathcal{R}=\mathcal{C}+V$. Take $v\in V$ and $l\in L$ such that $l-l_1\not\in V$. If $W_B(\tau_l(x_0+v))\neq0$ then, by Theorem \ref{thcora} $\tau_l(x_0+v)\in\tau_{l_1}x_0+V$. This implies that $\tau_l(x_0)-\tau_{l_1}x_0\in V$ so $\tau_{l}0-\tau_{l_1}0\in V$. With the hypothesis ``(H) modulo $V$'' we get $l-l_1\in V$, a contradiction. This proves the lemma.

\end{proof}
\begin{example}\label{ex1}

          To illustrate our method, we now give a natural but non-trivial
example $(R, B, L)$ in $\br^2$ for which $\mu_B$ may be seen to be a spectral
measure. In fact, we show that there is a choice for its spectrum
$\Lambda = \Lambda(\mu_B)$ which we compute with tools from Definition \textup{\ref{defredcond}},
Theorem \textup{\ref{thmain}}, and Lemma \textup{\ref{lem4_9}}.  Moreover, for the computation of the whole
spectrum $\Lambda$, the $W_B$-cycles do not suffice. \textup{(}There is one $W_B$ cycle, a
one-cycle, and it generates only part of $\Lambda$.\textup{)} Hence in this example,
the known theorems from earlier papers regarding spectrum do not suffice. To
further clarify the $W_B$-cycles in the example, we have graphed the two
attractors  $X_B$ and  $X_L$ in Figures \textup{\ref{figxb}} and \textup{\ref{figxl}}.

Take 
$$R:=\left[\begin{array}{cc}
4&0\\
1&4\end{array}\right],\quad B:=\left\{\vect{0}{0},\vect{0}{3}\vect10\vect13\right\}.$$
One can take 
$$L:=\left\{\vect00,\vect20,\vect02,\vect22\right\}.$$
One can check that the matrix in Definition \ref{defhada} is unitary so $(R,B,L)$ is a Hadamard triple.

\begin{figure}[h]
  \hfill
  \begin{minipage}[t]{.45\textwidth}
    \begin{center}  
\includegraphics[width=2in]{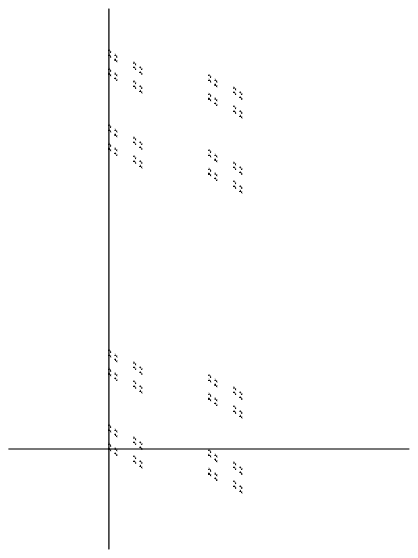}
      \caption{$X_B$}
      \label{figxb}
    \end{center}
  \end{minipage}
  \hfill
  \begin{minipage}[t]{.45\textwidth}
    \begin{center}  
\includegraphics[width=2in]{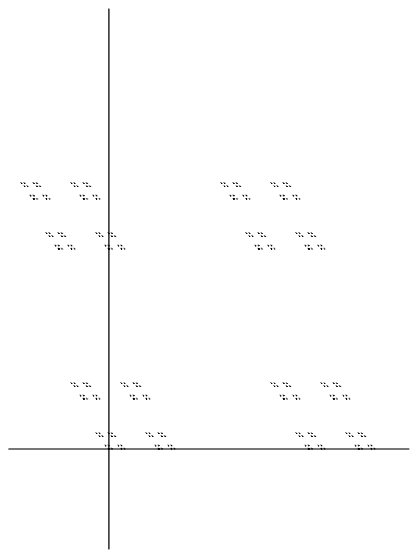}
      \caption{$X_L$}
      \label{figxl}
    \end{center}
  \end{minipage}
  \hfill
\end{figure}
\par
We look for $W_B$-cycles. We have
$$W_B(x,y)=\left|\frac{1}{4}(1+e^{2\pi i x}+e^{2\pi i 3y}+e^{2\pi i(x+3y)})\right|^2.$$
Then $W_B(x,y)=1$ iff $x\in\bz$ and $y\in\bz/3$ (all the terms in the sum must be equal to $1$). 
\par
If $(x_0,y_0)$ is a point of a $W_B$-cycle, then for some $(l_1,l_2)\in L$, $\tau_{(l_1,l_2)}(x_0,y_0)$ is also in the $W_B$-cycle, so $x_0,\frac{1}{4}(x_0+l_1)-\frac{1}{16}(y_0+l_2)\in\bz$ and $y_0,\frac{1}{4}(y_0+l_2)\in\bz/3$. Also, note that $(x_0,y_0)$ is in the attractor $X_L$ of the IFS $(\tau_l)_{l\in L}$, so $0\leq y_0\leq2/3$, and $-1/4\leq x_0\leq 2/3$. (This can be seen by checking that the rectangle $[-1/4,2/3]\times[0,2/3]$ is invariant for all $\tau_l$, $l\in L$.)
\par
Then, we can check these points and obtain that the only $W_B$-cycle is $(0,0)$, of length $1$, which corresponds to 
$\underline{\vect00}$.
\par
Now we look for the vector spaces $V$ that might appear in connection to the minimal invariant sets (see Theorem \ref{thcora}). Since these spaces are proper, and we have eliminated the case when $V=\{0\}$ by considering the $W_B$-cycles, it follows that $V$ must have dimension $1$ so it is generated by an eigenvector of $S=\left[\begin{array}{cc}4&1\\0&4\end{array}\right]$. Thus $V=\{(x,0)\,|\,x\in\br\}$.
\par
This subspace is reducible, with $r_1=0,r_2=1$, $\eta_{1,1}=0$, $\eta_{1,2}=3$, $\eta_{2,1}=0$, $\eta_{2,2}=3$, $s_1=0$, $s_2=2$, $\gamma_{1,1}=\gamma_{2,1}=0$, $\gamma_{1,2}=\gamma_{2,2}=2$. The measure $\mu_1$ on the first component corresponds to the IFS $\tau_0(x)=x/4$, $\tau_1(x)=(x+1)/4$. This corresponds to $R_1=4$, $B_1:=\{0,1\}$ and one can take $L_1:=\{0,2\}$ to get $(R_1,B_1,L_1)$ a Hadamard pair. The associated function is $W_{B_1}(x)=|\frac{1}{2}(1+e^{2\pi i x})|^2$. The only points where $W_{B_1}$ is $1$ are $x\in\bz$. Then one can see that the only $W_{B_1}$-cycle is $\{0\}$. Thus the spectrum of $\mu_1$ is 
$\{\sum_{k=0}^n4^ka_k\,|\,a_k\in\{0,2\}, n\in\bn\}$.
\par
We have to find the associated cycle $\mathcal{C}$. As in Lemma \ref{lemposinv}, we must have $W_B(\tau_l(x_0)+v)=0$ for elements $x_0$ in the cycle and some $l\in L$ and all $v\in V$. But this means that, for the second component $y'\in\br$ of $\tau_lx_0$, $1+e^{2\pi i x}+e^{2\pi i 3y'}+e^{2\pi i(x+3y')}=0$. This implies that $y'=(2k+1)/6$ for some $k\in\bz$. Moreover, we saw in Lemma \ref{lemwtildecycle} that the set of the second components of $\mathcal{C}$ must be a $\tilde W$ cycle. In our case $\tilde W(y)=\frac{1}{2}\left|1+e^{2\pi i3y}\right|^2$, and the IFS in case is $\{\tau_{s_i}\}=\{\tau_0,\tau_2\}$. The $\tilde W$-cycles are $\{0\}$ corresponding to $\underline{0}$, and $\{2/3\}$ corresponding to $\underline{2}$. Thus we obtain that the invariant sets obtained as translations of $V$ could be: $V$ and $\mathcal{R}:=2/3+V=\{(x,2/3)\,|\,x\in\br\}$. We can discard the first one because we see that $W_B(\tau_{(0,2)}(x,0))$ is not constant $0$. The set $2/3+V$ is indeed invariant, and we have $\tau_{(l_1,l_2)}(x,2/3)=0$ if $l_2=0$, and $\tau_{(l_1,l_2)}(x,2/3)\in 2/3+V$ if $l_2=2$.
\par
Next we want to compute the contribution of each of these invariant sets to the spectrum of $\mu_B$. 
\par
For the $W_B$-cycle $\{(0,0)\}$, of length $m=1$, we have as in Lemma \ref{lem4_9}, $$k_{\mathcal{C}}(\omega_1\wdots \omega_{k-1})=\omega_1+S\omega_2+\dots +S^{k-1}\omega_{k-1}$$ for all $\omega_1,\dots ,\omega_{k-1}\in L$.
By induction one can see that $S^n=\left[\begin{array}{cc}4^n&n4^{n-1}\\
0&4^n\end{array}\right]$. So the contribution from this $W_B$-cycle is
\[
\Lambda(0):=\left\{\left(\sum_{k=0}^n4^ka_k+g(b_0,\dots ,b_n),\sum_{k=0}^n4^kb_k\right)\biggm|_k,b_k\in\{0,2\}\right\},
\]
where $g(b_0,\dots ,b_n)=\sum_{k=0}^nk4^{k-1}b_k$.
\par
For the invariant set $\mathcal{R}=\{(x,2/3)\,|\,x\in\br\}$, we have as in Lemma \ref{lem4_9}, with $\omega_1,\dots ,\omega_{k-1}\in\{0,2\}$,
$k_{2/3}(\omega_1,\dots ,\omega_{k-1})=\omega_1+4\omega_2+\dots +4^{k-1}\omega_{k-1}-4^k\frac23$, or writing $2/3=2/4+2/4^2+\dots +2/4^{k}+2/4^{k+1}+\cdots $, we obtain
$$k_{2/3}(\omega_1,\dots ,\omega_{k-1})=\sum_{i=0}^{k-1}a_i4^{i}-\frac23,$$ with $a_k\in\{0,-2\}$.
\par
As in Proposition \ref{propnr} and its proof, using the spectrum of $\mu_1$, the contribution to the spectrum is 
\[
\Lambda(2/3):=\left\{\left(\sum_{k=0}^n4^ka_k,-\frac23-\sum_{k=0}^m4^kb_k\right)\biggm|a_k,b_k\in\{0,2\},n,m\in\bn\right\}.
\]
\par

Finally, the spectrum of $\mu_B$ is $\Lambda_B:=\Lambda(0)\cup\Lambda(2/3)$.
\par
Note also, that we can use the decomposition given in Lemma \ref{lem4_4}. The measure $\mu_1$ is the invariant measure for the IFS: $\tau_{0}(x)=x/4$, $\tau_1(x)=(x+1)/4$. For all $x\in\br$, the measure $\mu_x^2=:\mu_2$ is the invariant measure for the IFS $\tau_{0}(x)=x/4$, $\tau_3(x)=(x+3)/4$. Both $\mu_1$ and $\mu_2$ are spectral measures (one can use $L=\{0,2\}$ for both of them). We saw that the spectrum of $\mu_1$ is $\Lambda_1:=\{\sum_{k=0}^n4^ka_k\,|\,a_k\in\{0,2\},n\in\bn\}$. The IFS $(\tau_{0},\tau_{3})$ has two $W_{B_2}$-cycles: $\{0\}$ and $\{2/3\}$, so, after a computation we get that the spectrum of $\mu_2$ will be $\Lambda_2:=\Lambda_1\cup(-\frac23-\Lambda_1)$. 
\par
Using the decomposition of Lemma \ref{lem4_4} we obtain that a spectrum for $\mu_B$ is $\Lambda_1\times\Lambda_2$. It is interesting to see that this is a different spectrum than the one computed before $\Lambda_B$. 
\end{example}

\begin{remark}
Since in Example \ref{ex1}, the $W_B$-cycles are not sufficient to describe all invariant sets, the results from \cite{JoPe98,Str00,LaWa02,DuJo05} do not apply here; they give only part of the spectrum, namely the contribution of the $W_B$-cycle $\{0\}$. 
\end{remark}

\begin{example}
Take now $B$ to be a complete set of representatives for $\bz/R\bz^d$. So $N=|\det R|$. To get a Hadamard triple, one can take $L$ to be any complete set of representatives for $\bz^d/S\bz^d$, because the matrix $\frac{1}{\sqrt{N}}(e^{2\pi i b\cdot l})_{b,l}$ will then be the matrix of the Fourier transform on the finite group $\bz^d/R\bz^d$, hence unitary. 
\par
The following proposition is folklore for affine IFSs; see, e.g., \cite{CHR97, JoPe94, JoPe96, LaWa96}.
\begin{proposition}
Suppose the vectors in $B$ form a complete set of coset representatives for the finite group $\bz^d/R\bz^d$.
Then the following conclusions hold:
\begin{enumerate}
\item[(a)] The attractor $X_B$ has non-empty interior relative to the metric from $\br^d$.
\item[(b)] The Borel probability measure $\mu_B$ is of the form $\mu_B =\frac{1}{p}$(Lebesgue measure in $\br^d$ restricted to $X_B$), where $p$ is an integer.
\item[(c)] Moreover, $p = 1$ if and only if the attractor $X_B$ tiles $\br^d$ by translations with vectors in the standard lattice $\bz^d$; where by {\it tiling} we mean that the union of translates $\{X_B + k\, |\,  k\in\bz^d\}$ cover $\br^d$ up to measure zero, and where different translates can overlap at most on sets of measure zero.
\item[(d)] In general, there is a lattice $\Gamma$ contained in $\bz^d$ such that $X_B$ tiles $\br^d$ with $\Gamma$; and the group index $[\bz^d : \Gamma]$ coincides with the number $p$.
\end{enumerate}

\end{proposition}
Using Fuglede's theorem \cite{Fug74} it follows that $\mu_B$ is a spectral measure, with spectrum the dual lattice of $\Gamma$. (Fuglede's theorem \cite{Fug74} characterizes measurable subsets $X$ in $\br^d$ which are fundamental domains for some fixed rank-$d$ lattice $L$. First note that such subsets have positive and finite Lebesgue measure,  $\mu =$ the $d$-dimensional Lebesgue measure. For measurable fundamental domains, Fuglede showed that $L^2(X, \mu)$ has $\{e_\lambda\, |\,  \lambda \mbox{ in the dual lattice to } L\}$ as ONB, i.e., that the dual lattice is a set of Fourier frequencies. More importantly, he proved the converse as well: If $L^2(X, \mu)$ for some measureable subset of $\br^d$ is given to have an ONB consisting of a lattice of Fourier frequencies, then $X$ must be a fundamental domain for the corresponding dual lattice. Furthermore, he and the authors of \cite{Ped87,JoPe92} also considered extensions of this theorem to sets of Fourier frequencies that are finite unions of lattice points. We should add that there is a much more general Fuglede problem which was shown recently \cite{Tao04} by Tao to be negative.) 

The relation between the lattice $\Gamma$ and the invariant sets will be the subject of another paper.
\end{example}

\addvspace{\medskipamount}\noindent\textbf{Notes on the literature.}
While there is, starting with \cite{Hut81} and \cite{BEHL86}, a substantial literature
of papers treating various geometric features of iterated function systems
(IFS), the use of Fourier duality is of a more recent vintage. The idea of
using substitutions together with duality was perhaps initiated in \cite{JoPe92};
see also \cite{Mas94}. However, the use of substitutions in dynamics is more
general than the context of IFSs; see, for example, \cite{LiMa95}. We further want
to call attention to a new preprint \cite{Fre06} which combines the substitution
principle with duality in a different but related manner. The use of duality
in \cite{Fre06} serves to prove that the class of affine IFSs arises as model
sets. It is further interesting to note (e.g., \cite{Bar01}) that these fractals
have found use in data analysis.
\par
In the definition of reducible subspaces we added a certain non-overlapping condition for the measure $\mu_1$. This condition, which might be automatically satisfied for our affine IFSs, is part of a more general problem:
\begin{problem}
Give geometric conditions for a fixed $(X, \tau_i)$ which guarantee that the distinct sets $\tau_i(X)$ overlap at most on subsets of $\mu$-measure zero.
\end{problem}
For related but different questions, the reader can consult \cite{Sch94,LaWa93,LaRa03,HLR03}.

\begin{acknowledgements}
We gratefully acknowledge discussions
with professors Yang Wang and Ka-Sing Lau. In addition, this work was supported by
a grant from the National Science Foundation (NSF-USA), DMS 0457491. The authors thank Brian Treadway for expert help with tex problems and with graphics. The co-authors are very grateful to the referee for his/her careful work on
our manuscript and his/her thoughtful suggestions.
We have followed them all, and we are grateful to him for sharing them with
us.

\end{acknowledgements}


\end{document}